\documentclass[]{amsart} 
\usepackage[]{amsmath, amsthm, amsfonts}
\usepackage[all]{xy} 
\usepackage{hyperref} 

\newcommand {\C} {{\mathbb C}}

 \newcommand {\Q} {{\mathbb Q}} 
\newcommand {\PP} {{\mathbb P}}

 \newtheorem{thm}{Theorem}
 \newtheorem{cor}{Corollary}
 \newtheorem{lemma}{Lemma}

\begin{document}

 \title{The Hodge conjecture for rationally connected fivefolds }

 \author{ Donu Arapura} \thanks{Author partially supported by the NSF}
 \address{Department of Mathematics\\
   Purdue University\\
   West Lafayette, IN 47907\\
   U.S.A.}
 
 \maketitle

{\small It has been pointed out to me, by A. Collino and H. Esnault,
that the road taken here has been well travelled. The papers of 
Bloch-Srivinas \cite{bs},
Esnault-Levine \cite{el}, Fakhruddin-Rajan \cite{fr} --
especially the first  -- contain a
number of related ideas and applications, if not precisely 
the results proven here. That said, I will let the original text
stand as it was written. Since the flatness of its prose
may not quite have hidden the pleasure of writting it, I hope that
this note can at least serve as introduction to this beautiful
technique.}

The purpose of this note is to prove the Hodge conjecture for
a five dimensional smooth projective complex variety $X$ 
possessing sufficiently many rational
curves.  To make this precise, recall \cite{kollar} that
there is a rational map $f:X\dashrightarrow Y$ called the maximal
rationally connected fibration. The fibers of $f$ are rationally connected,
and $\dim Y$ is minimal among all such maps.
We show that the Hodge conjecture holds provided that $\dim Y\le 3$.
In particular, it holds if $X$ is already rationally connected as would
be the case if it were Fano.

As to why we state the result only in dimension 5, we should
recall  that  the Hodge conjecture is  only an issue in dimensions 
greater than 3, but  for fourfolds covered by rational curves the result  
is well known and elementary.

\begin{lemma}[Conte-Murre]\label{lemma:conte}
  The Hodge conjecture holds for a smooth projective uniruled fourfold
over $\C$.
\end{lemma}

The idea is that  since a uniruled fourfold is dominated by a blow up
of a product of $\PP^1$ with a threefold $Y$, its interesting cohomology 
``comes from'' the lower dimensional varieties
$Y$ and the centers of the blow ups. 
We will give a more precise explanation later.
 A similar argument would be sufficient 
to  prove  the main theorem when  $f:X\dashrightarrow Y$
is fiberwise unirational i.e. if the map is dominated by a 
blow up of $\PP^{5-\dim Y}\times Y$. 
In general, however, a more sophisticated argument is  needed.
The key trick  is the following lemma of Bloch that we learned  from
the papers of Esnault \cite{esnault} and Kim \cite{kim}.

 \begin{lemma}[Bloch]
   Let $X$ be a smooth rationally connected variety over an algebraically
closed field $K$,
then a positive multiple of the diagonal $\Delta\subset X\times X$
is rationally equivalent to a sum $\xi\times X+ \Gamma$,
where $\xi\in Z_0(X)$ is a zero cycle,  
and $\Gamma$ is supported on $X\times Z$
for some proper closed subset $Z\subset X$.
 \end{lemma}

 The original statement is a bit  stronger \cite{bloch, esnault}, 
but we won't need it. In order to keep our
presentation reasonably self contained, we give a geometric proof of the
 result in its present form.

 \begin{proof}
Rational connectedness means that any two general points can be
joined by an irreducible rational curve.
More precisely \cite[thm II 2.8, IV def 3.2]{kollar},
there is a morphism $R:\PP^1\times F\to X$,
such that
$$R^{(2)}:\PP^1\times \PP^1\times F\to X\times X,\> 
R^{(2)}(t_1,t_2,f) = (R(t_1,f), R(t_2,f))$$
is dominant.
Fix a general point $x_0\in X$. Consider the
maps $F\to X\times X$ given by the composition of $0\times \infty\times id_F$
and $R$, and $X\to X\times X$ given by $x_0\times id_X$.
Let $U= F\times_{X\times X} X$ be the fiber product.
The restriction of $R$ gives a morphism $r:\PP^1\times U\to X$, 
with $r(0,u)=x_0$ such that
$g:U\to X$ given by $g(u) = r(\infty,u)$
is dominant. After replacing
$U$ by the normalization of $K(X)$ in $K(U)$, we can
assume that $\dim U =   \dim X$.
Thus we have a commutative diagram
 $$
\xymatrix{
 \PP^1\times U\ar[r]^{q}\ar[d]^{\pi_2} & X\times X\ar[d]^{\pi_2} \\ 
 U\ar[r]^{g} & X
}
$$
where $\pi_2$ are projections to the second factor,
and $q(t,u) = (r(t,u), g(u))$.  We see that  $q(\infty\times U) =
\Delta\cap (X \times g(U))$ and $q(0\times U)= x_0\times g(U)$ as sets.
Therefore the difference  $deg(g)(\Delta-x_0\times X)
\in CH_*(X\times X)$ is supported on  $X\times \overline{(X-g(U))}$. 
 \end{proof}

\begin{cor}\label{cor:bloch}
 Suppose that $f:X\to Y$ is a morphism of varieties such that the fibers
of $f$ are rationally connected. Then
 there exists a nonempty open subset $V\subset Y$, a relative
zero cycle $\xi\in Z_0(X\times_Y V/V)$,
 a proper closed subset $Z\subset X_V=X\times_Y V$, and a cycle
 $\Gamma$
supported
on $X\times_Y Z$ such that a multiple of the relative diagonal
$\Delta\subset X_V\times_V X_V$ is rationally equivalent to
 $\xi\times X+ \Gamma$ after restriction to the (closed) fibers of $X_V\to V$.
 \end{cor}

 \begin{proof}
The geometric general fiber $X_{\bar \eta}=X\times_Y\, Spec
\overline{K(Y)}$ is rationally connected by (the proof
of) \cite[thm IV 3.11]{kollar}. 
The previous lemma applied to the $X_{\bar \eta}$, produces
cycles, $\xi'$, $Z'$, $\Gamma'\subset X\times Z'$,  on 
$X_{\bar  \eta}\times X_{\bar\eta}$,
such that $\xi'\times X_{\bar \eta}+ \Gamma'$
is rationally equivalent to  the diagonal $\Delta'$ times a positive
integer $N$.
These cycles and the data defining the rational equivalence are
defined over a finite extension of $K(Y)$. Therefore they can
be spread out to relative cycles $\xi''$, $Z''$, $\Gamma''$ and data
on  $X_{\tilde V}$, for some generically finite map $\tilde V\to Y$. 
We then have  $N$ times the relative diagonal of $X_{\tilde V}$ 
is rationally equivalent
to $\xi''\times X_{\bar \eta}+ \Gamma''$ on the fibers of $X_{\tilde V}$.
Now let $\xi, Z,\Gamma$ be  the pushforward of $\xi'', Z'',\Gamma''$ onto 
$X_V$ where $V = im(\tilde V)$. Then $(deg \tilde V/V)N\Delta$
will be rationally equivalent to $\xi\times X+ \Gamma$ on the fibers,
and the support of $\Gamma$ will lie in $X\times_Y Z$.
 \end{proof}

From now on, we assume our varieties are complex. Cohomology
is singular cohomology with rational coefficients. These groups carry
mixed Hodge structures, but we will omit Tate twists, unless they
seem essential.
Given a smooth projective variety $X$ and closed subset $Z$, 
we want to be able to reduce the verification of the Hodge conjecture 
on $X$ into a separate verifications  on $Z$ and $X-Z$.
Here we use Jannsen's  formulation \cite{jannsen} 
for  the Hodge conjecture for singular
or nonproper varieties, which says that the space of rational
$(-p,-p)$ cycles in the pure part of Borel-Moore homology 
$Gr_{-2p}^WH_{2p}(X)$ consists of  algebraic cycles. In other words, 
this space is  spanned
by fundamental classes of $p$-dimensional subvarieties of $X$.
When $X$ is smooth, we can formulate this
in compactly supported cohomology  thanks to the duality isomorphism 
$$H_c^{2(\dim X-p)}(X,\Q)\cong H_{2p}(X,\Q)(-\dim X).$$
It is worth noting that this form of the
Hodge conjecture for general $X$ 
 would follow from the usual Hodge conjecture on a desingularization
of a compactification of $X$ \cite[thm 7.9]{jannsen}.\\ \\

\begin{lemma}\label{lemma:HCforOpen}
  \begin{enumerate}
\item[]
  \item  If $X$ is an open subset of a smooth
projective variety $\bar X$, then the Hodge conjecture
holds for  $H^{2p}_c(X)$, i.e. 
rational $(p,p)$ cycles
in it are all algebraic if the conjecture
holds for $H^{2p}(\bar X)$.
\item 
If $X$ is projective with a desingularization $\tilde X\to X$,
then the Hodge conjecture holds for $H_{2p}(X)$ i.e. rational $(-p,-p)$ cycles
in it are all algebraic if the same statement
holds for $H_{2p}(\tilde X)$ i.e. if the usual form of the conjecture
holds for $H^{2(\dim X-p)}(\tilde X)$.

\item Suppose that $X$ is a projective variety with a closed subset 
$Z\subset X$.
If the Hodge conjecture holds  for $H_{2p}(Z)$ and $H_{2p}(X-Z)$, 
then it holds for   $H_{2p}(X)$. 
  \end{enumerate}
\end{lemma}

\begin{proof}
The proof of the first two statements are essentially contained in 
the arguments on \cite[pp. 113-114]{jannsen}. In summary, for (1),
we use the injectivity of $Gr_{2p}^WH_c^{2p}(X)\to H^{2p}(\bar X)$
to reduce the Hodge conjecture for $H_c^{2p}(X)$ to the corresponding 
statement $H^{2p}(\bar X)$.
For (2), we use the surjectivity of
$H_{2p}(\tilde X)\to Gr_{-2p}^W H_{2p}(X)$ to achieve a similar reduction.

For the third statement, we use the exact sequence of mixed Hodge structures
$$H_{2p}(X-Z)\to H_{2p}(X)\to H_{2p}(Z)$$
to obtain an exact sequence of pure Hodge structures
$$Gr_{-2p}^WH_{2p}(Z)\to Gr_{-2p}^W H_{2p}(X)\to Gr_{-2p}^W H_{2p}(X-Z),$$
and then a similar sequence for the spaces of Hodge i.e.
rational $(-p,-p)$-cycles.
Given a Hodge cycle $\alpha$ in $Gr_{-2p}^WH_{2p}(X)$, its image $\beta$ in 
$Gr_{-2p}^WH_{2p}(X-Z)$ is algebraic by hypothesis. Thus
$\beta=\sum_i\, n_i [V_i]$.
 Taking closures of the components allows us to lift $\beta$
to an algebraic cycle $\bar \beta=\sum_i\, n_i[\bar V_i]$ on $X$. 
Then the difference
$\beta-\bar \beta$ would  be represented by a Hodge and hence algebraic cycle
supported on $Z$.
\end{proof}

As a prelude to the proof of the theorem, let us sketch a proof
of lemma~\ref{lemma:conte}. 

\begin{proof}[Proof of lemma 1]
Given a uniruled fourfold $X$,
it is dominated by 
a blow up $\tilde X$ of $\PP^1\times Y$ with $\dim Y=3$.
Since any Hodge cycle can be lifted to $\tilde X$, it
suffices to prove the conjecture for this. Now $\tilde X$
can be decomposed into a disjoint union of strata of
the form $\C^i\times S$ with $\dim S\le 3$.
Part 3 of the previous lemma shows that it suffices to check the
Hodge conjecture for these strata. By part 1 of the lemma, it
suffices to do this for compactifications
$\PP^i\times \bar S$. But this follows immediately
from K\"unneth's formula.
\end{proof}

\begin{thm}
Let $X$ be a smooth projective  five dimensional variety over $\C$, such
that the base of the maximal rationally connected fibration is
at most three dimensional. Then the  Hodge conjecture holds for $X$.
 \end{thm}

 \begin{proof}
 It is enough to check the Hodge conjecture for $H^n(X)$ for $n=4$,
since the other cases would follow from the Lefschetz $(1,1)$ theorem
when $n=2$,
or by previous cases and  Hard Lefschetz when $n\ge 6$.
Let $f:X\dashrightarrow Y$ be maximal rationally connected fibration.
Since the validity of the Hodge conjecture for $X$ is implied by 
its  validity on any smooth blow up $\tilde X\to X$, we
can assume that $f$ is a morphism. Note that all the fibers of $f$
are rationally connected  by  \cite[thm IV 3.11]{kollar}.

We now apply  corollary \ref{cor:bloch} to conclude
that the relative diagonal $\Delta\subset X_V\times_V X_V$ decomposes as
$\xi\times X+ \Gamma$ in the Chow groups $CH_*(X_y\times X_y)\otimes \Q$
with $y\in V$. 
We let  $D\subset X$  be a divisor containing $Z$, and let $\tilde D\to
D$ be a desingularization. Note that we can assume that $f_V:X_V\to V$ 
and $g:\tilde D_V\to V$ are smooth, after shrinking $V$.
Let $W=Y-V$.
By lemma \ref{lemma:HCforOpen}, it is enough to check the Hodge conjecture
for $H_c^4(X_V)\cong H_{6}(X_V)$ and $H_{6}(X_{W})$ separately. For 
$X_W$, note that it is a uniruled variety of dimension at most 4.
The same goes for any desingularization of it. Therefore $X_W$
satisfies the Hodge conjecture in all degrees by lemma \ref{lemma:conte}.

Now we turn to $X_V$.
The identity on $H^i(X_y,\Q)$, with $y\in V$, is given by the action of the
restriction of the correspondence
$\Delta$. This decomposes into a sum of the action of the restriction
of $ \xi\times X$,
which vanishes for $i>0$, and the action of the restriction of $\Gamma$.  
Note that the image under $\Gamma|_{X_y}$ lies in the
kernel of $H^i(X_y,\Q)\to H^i(X_y-D_y,\Q)$, and this coincides with
the image of the Gysin
map $H^{i-2}(\tilde D_y,\Q)(-1)\to H^i(X_y,\Q)$ by \cite[Prop 8.2.8]{deligne}.
It follows that the map of local systems  $R^{i-2}g_*\Q\to
R^if_{V*}\Q$ is surjective; it is necessarily split surjective by 
Deligne's theorem on semisimplicity of monodromy \cite[thm
4.2.6]{deligne}. Therefore
$$H_c^j(V, R^{i-2}g_*\Q) \to H_c^j(V, R^{i}f_{V*}\Q)$$
is surjective for $i>0$. This together with the fact that
 the Leray spectral sequence degenerates \cite{deligne-L} implies
that $H_c^4(X_V,\Q)$ is spanned by the sum of the images of 
$H_c^2(\tilde D_V,\Q)(-1)$ and $H_c^4(V,\Q)$ as  mixed Hodge structures. 
In particular, it suffices to check the Hodge conjecture for 
$H_c^2(\tilde D_V,\Q)$ and $H_c^4(V,\Q)$. By lemma \ref{lemma:HCforOpen},
we can check the corresponding statements after 
 replacing $\tilde D_V$ and $V$ by their compactifications $\tilde D$ and
$Y$. The Hodge conjecture holds for $H^2(\tilde D,\Q)$ by 
the Lefschetz $(1,1)$ theorem, and it holds for $Y$ since
$\dim Y\le 3$. Now we are done.
\end{proof}
 
Recall that smooth variety is Fano if its anticanonical bundle $\omega^{-1}$
is ample. Such varieties are rationally connected by the
work of Campana, Koll\'ar, Mori and Miyaoka \cite[cor. V 2.14]{kollar}. Thus:

\begin{cor}
   The Hodge conjecture holds for a Fano fivefold.
 \end{cor}

 \begin{cor}
  The Hodge conjecture holds for a smooth hypersurface
in $\PP^n\times \PP^{6-n}$ of bidegree $(a,b)$ (or degree
$a$ when $n=6$) with 
$3\le n\le 6$ and $a\le n$.
 \end{cor}

 \begin{proof}
   Projection of the hypersurface to $\PP^{6-n}$ yields a family of 
 Fano varieties on it of dimension  two or more. This forces the base of the
maximal rationally connected family to be less than or equal to $3$.
 \end{proof}

By similar arguments, it follows that a rationally connected $n$-fold
is motivated by an $(n-2)$-fold in the sense of \cite{arapura}. In
particular, the Lefschetz standard conjecture holds for a rationally
connected fourfold.

\end{document}